\documentclass[hidelinks,onefignum,onetabnum,final]{siamart220329}

\usepackage{lipsum}
\usepackage{amsfonts}
\usepackage{bm}
\usepackage{graphicx}
\usepackage{epstopdf}
\usepackage{algorithmic}
\ifpdf
  \DeclareGraphicsExtensions{.eps,.pdf,.png,.jpg}
\else
  \DeclareGraphicsExtensions{.eps}
\fi

\newsiamremark{remark}{Remark}
\newsiamremark{hypothesis}{Hypothesis}
\crefname{hypothesis}{Hypothesis}{Hypotheses}
\newsiamthm{claim}{Claim}

\headers{RBF-FD boundary discretization for elliptic PDEs}{M. E. Nielsen and B. Fornberg}

\title{High-order numerical method for solving elliptic partial differential equations on unfitted node sets}

  \author{Morten E. Nielsen\thanks{DHI, Hørsholm, 2970, Denmark }
\and Bengt Fornberg\thanks{Department of Applied Mathematics, University of Colorado, Boulder,  CO 80309, USA}}

\usepackage{amsopn}

\ifpdf
\hypersetup{
  pdftitle={High-order numerical method for solving elliptic partial differential equations on unfitted node sets},
  pdfauthor={M. E. Nielsen and B. Fornberg}
}
\fi

\begin{document}

\maketitle

\begin{abstract}
 In this paper, we present how high-order accurate solutions to elliptic partial differential equations can be achieved in arbitrary spatial domains using radial basis function-generated finite differences (RBF-FD) on unfitted node sets (i.e., not adjusted to the domain boundary). In this novel method, we only collocate on nodes interior to the domain boundary and enforce boundary conditions as constraints by means of Lagrange multipliers. This combination enables full geometric flexibility near boundaries without compromising the high-order accuracy of the RBF-FD method. The high-order accuracy and robustness of two formulations of this approach are illustrated by numerical experiments.
\end{abstract}

\begin{keywords}
High-order, RBF-FD, Meshfree, PDEs, Lagrange multipliers
\end{keywords}

\begin{MSCcodes}
65N22, 35J05  
\end{MSCcodes}

\section{Introduction}

Partial differential equations (PDEs) form key parts of virtually all mathematical models describing physical systems 
and processes, with a large number of numerical approaches available. This paper focuses on accurate treatment of 
boundaries of non-trivial (and possibly time dependent) shapes. A first key numerical choice is whether to use a grid-based or 
a meshfree PDE discretization:

\vspace{0.25cm}

\begin{itemize}
    \item Grid-based: The main options are overlapped grids (\cite{chesshire1990composite,henshaw2006moving}) and immersed boundary techniques (\cite{zapata2023compact,li2006immersed,mittal2005immersed,peskin1977numerical,peskin2002immersed})
    \item Meshfree: The main option up to now has been to use quasi-uniform node sets that align with boundaries (e.g., radial basis function methods (\cite{FF2015_RBFPrimer,fornberg2015solving})). The present novelty is to use
‘unfitted’ node sets that do not impose the constraint of node sets being aligned with boundaries.
\end{itemize}

\vspace{0.25cm}
This last option, meshfree quasi-uniform node sets not fitted to align with boundaries is attractive because, it

\vspace{0.25cm}

\begin{itemize}
    \item Simplifies the generation of node sets,
    \item Offers the most flexibility in permitting local refinement purely based on accuracy considerations,
    \item Simplifies the handling of moving boundaries (as no corresponding movement of nodes is called for).
\end{itemize}

\vspace{0.25cm}

High-accuracy implementations with unfitted quasi-uniform node sets have so far not been explored. We demonstrate 
here that combining RBF-FD (radial basis function-generated finite differences) with a Lagrange multiplier approach 
at boundaries can achieve very high convergence rates (here tested up to $\mathcal{O}(h^{10})$, where $h$ is an average node spacing, or equivalently in 2-D, $\mathcal{O}(1/N^{5})$ where $N$ is the total number of nodes in a discretization. We here 
introduce and test this new approach for Poisson problems in two- and three-dimensional space with different boundary conditions. Computationally optimized implementations and 
other types of PDEs are left to be considered for future research.

\begin{figure}[h!]
	\centering
		\includegraphics[scale=0.2]{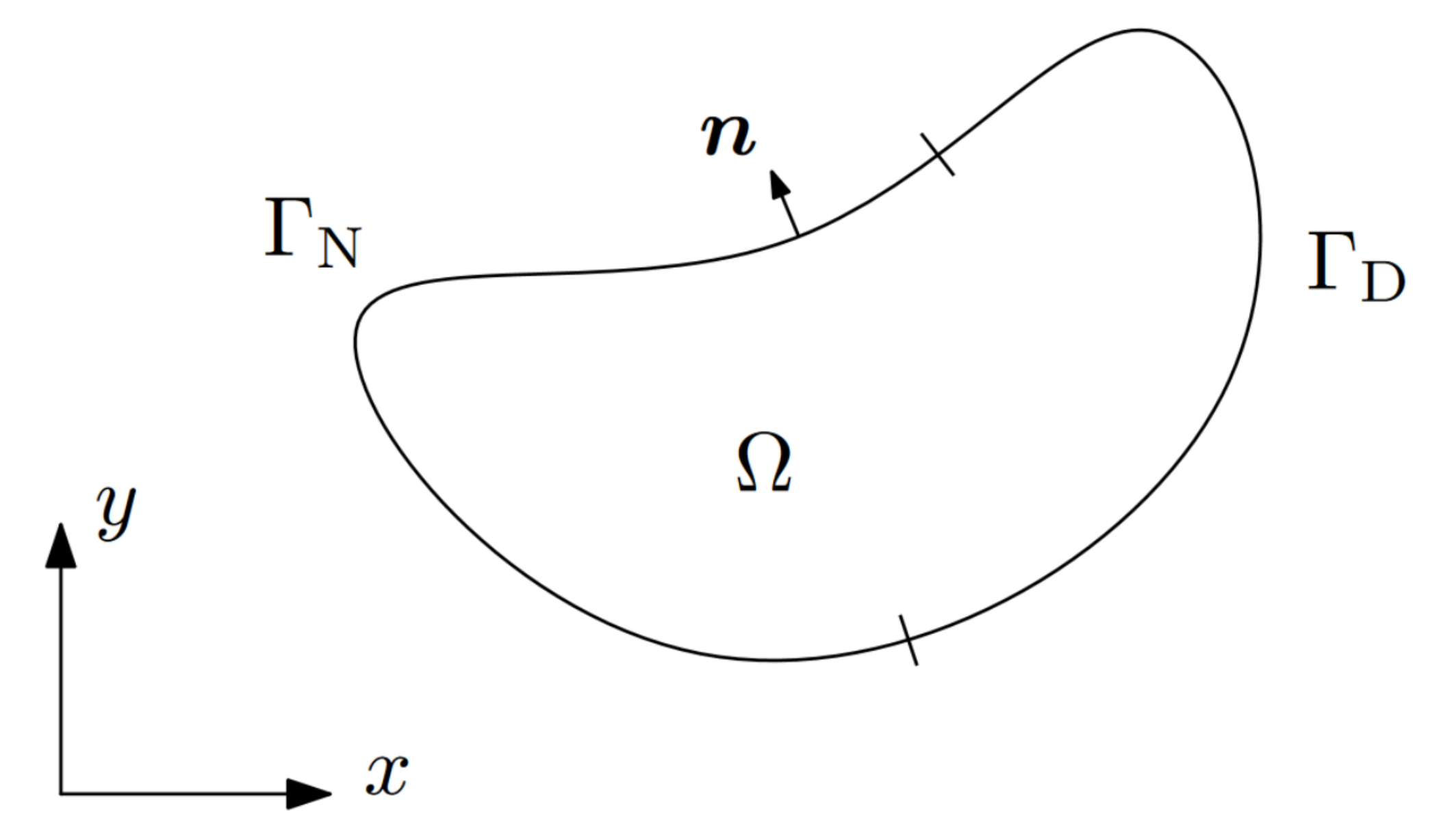}
	\caption{Domain with boundary definitions for $d = 2$.}
	\label{fig:domain}
\end{figure}

In this work elliptic PDEs are considered and, for the sake of simplicity, only different variations of the Poisson problem are tested. Let $\Omega \in \mathbb{R}^d$ be a $d$-dimensional domain bounded by $\Gamma \in \mathbb{R}^{d-1}$, in which we seek a solution to the Poisson problem,

\begin{equation}\label{eq:peq}
\begin{aligned}
\nabla^2 u &= f(\bm{x}),\hspace{0.25cm} \bm{x} \in \Omega, \\
\mathcal{B}_{\text{D}} (u) = u &= g(\bm{x}), \hspace{0.25cm} \bm{x} \in \Gamma_{\text{D}}, \\
\mathcal{B}_{\text{N}} (u) = \frac{\partial u}{\partial n} &= h(\bm{x}),\hspace{0.25cm} \bm{x} \in \Gamma_{\text{N}},
\end{aligned}
\end{equation}

\noindent where $\nabla^2 \in \mathbb{R}$ is the continuous Laplace operator, $u = u(\bm{x}) \in \mathbb{R}$ is the solution in $\Omega$, while $\Gamma_{\text{D}}$ specifies the Dirichlet boundary and $\bm{n}$ is the outward pointing normal vector defined on the Neumann boundary $\Gamma_{\text{N}}$ as depicted in figure \ref{fig:domain}.

In the following section, section 2, we introduce the traditional collocation-based RBF-FD method along with the novel Lagrange multiplier-based RBF-FD methods. Section 3 introduces several test problems and gives numerical comparison results, which is followed by section \ref{sec:concl} with some conclusions.

\section{Methodology} \label{sec:method}

\begin{figure}[h!]
	\centering
		\includegraphics[scale=0.75]{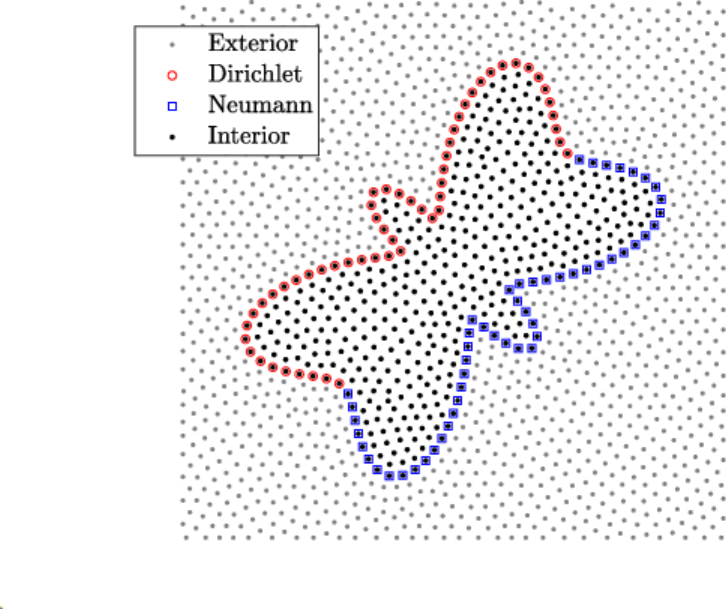} \vspace{-1cm}
	\caption{Boundary-fitted node set including definitions of exterior (gray dots), Dirichlet (red circles), Neumann (blue squares) and interior (black dots) boundary nodes. This node set combines quasi-uniformity with nodes being boundary-fitted.}
	\label{fig:fitnodes}
\end{figure}

The objective of this paper is to solve problems of the form (\ref{eq:peq}) with high-order accuracy while preserving as much of the meshfree nature of RBF-FD as possible. In the traditional RBF-FD collocation (RBF-FD-C) method, the computational domain is usually discretized by collocating $\text{N}_\text{I}$ nodes in $\Omega$, $\text{N}_\text{D}$ nodes on $\Gamma_{\text{D}}$, and $\text{N}_\text{N}$ nodes on $\Gamma_{\text{N}}$ according to the definitions illustrated in figure \ref{fig:domain}. An example of a boundary-fitted, quasi-uniform node set\footnote{Non-grid-based, locally uniform node set with possibly globally varying node density.} is illustrated in figure \ref{fig:fitnodes}. For the sake of completeness, the RBF-generated finite differences are first described briefly in section \ref{sec:rbffd} before being utilized.

\subsection{Radial basis function-generated finite differences} \label{sec:rbffd}

\noindent 
In the RBF-FD approach, the first step is to create FD-type stencils, one for each derivative and for each of the $N$ nodes within the domain. Including the stencil ‘evaluation’ point (at which it is to be accurate), we determine weights to 
use at its $n$ nearest neighbors (with $n \ll N$) as illustrated in figure \ref{fig:nset_test1}. These RBF-FD weights are obtained by creating an RBF-based interpolant across these $n$ nodes and differentiating this analytically. Weights at the nodes are thus obtained similarly 
to traditional 1-D FD weights when instead using polynomial interpolants \cite{FF2015_RBFPrimer,FLYER201639,bayona2017role}.

\noindent We introduce the polyharmonic radial basis functions $\phi_i(r)= \phi \left(\|\bm{x} - \bm{x}_i\|_2 \right) = \|\bm{x} - \bm{x}_i\|_2^{2k+1}$ and multivariate monomials $p_j(\bm{x})$ to approximate the exact solution to (\ref{eq:peq}) as,

\begin{equation} \label{eq:ansatz}
    u(\bm{x}) \approx u_h(\bm{x}) =  \sum_{i=1}^{n} \kappa_i \phi \left(\|\bm{x} - \bm{x}_i\|_2 \right) + 
    \sum_{j=1}^{\ell} \gamma_j p_j(\bm{x})
\end{equation}

\noindent  which we require to match at $n$ nodes, i.e. spatial points $\{\bm{x}_i\}_{i=1}^{n}$,

\begin{equation}
u(\bm{x}_i) = u_h(\bm{x}_i), \hspace{0.25cm} \text{for} \hspace{0.25cm} i=1, 2, ..., n,
\end{equation}

\noindent  while enforcing the additional constraints,

\begin{equation}
\sum_{i=1}^{n} \kappa_i p_j(\bm{x}_i) = 0, \hspace{0.25cm} \text{for} \hspace{0.25cm} j=1, 2, ..., \ell.
\end{equation}

 \noindent Here $\ell$ is the number of monomial terms in a multivariate polynomial of degree $m$. The above equations can be arranged in a linear system of equations,

\begin{equation} 
\tilde{A}     \begin{bmatrix}
        \bm{\kappa} \\
        \bm{\gamma}
    \end{bmatrix}
    =
    \begin{bmatrix}
        A & P\\
        P^T & 0
    \end{bmatrix}
    \begin{bmatrix}
        \bm{\kappa} \\
        \bm{\gamma}
    \end{bmatrix}
    =
    \begin{bmatrix}
        \bm{u} \\
        \bm{0}
    \end{bmatrix},
\end{equation}

\noindent where $\bm{\kappa}, \bm{u} \in \mathbb{R}^n$, $\bm{\gamma} \in \mathbb{R}^\ell$, $A_{ij} = \phi \left(\|\bm{x}_i - \bm{x}_j\|_2 \right)$ is an entry of the RBF collocation matrix $A \in \mathbb{R}^{n \times n}$ and $P_{ij}= p_j(\bm{x}_i)$ is an entry of the supplementary polynomial matrix $P \in \mathbb{R}^{n\times \ell}$. The linear operation $\mathcal{L}$ can be approximated at an evaluation point $\bm{x}_e$ as,

\begin{equation}
\mathcal{L} u(\bm{x})|_{\bm{x}_e} \approx  \mathcal{L} u_h(\bm{x})|_{\bm{x}_e} = \sum_{i=1}^{n} \kappa_i \mathcal{L}\phi \left(\|\bm{x} - \bm{x}_i\|_2 \right)|_{\bm{x}_e} + 
    \sum_{j=1}^{\ell} \gamma_j \mathcal{L} p_j(\bm{x})|_{\bm{x}_e},
\end{equation}

\noindent which again can be arranged in matrix-vector format,

\begin{equation}
   \mathcal{L} u(\bm{x})|_{\bm{x}_e}  \approx  
   \left[  \bm{a}^T \hspace{0.2cm} \bm{b}^T \right]     
   \begin{bmatrix}
        \bm{\kappa} \\
        \bm{\gamma}
    \end{bmatrix}
    =
    \left[  \bm{a}^T \hspace{0.2cm} \bm{b}^T  \right]  \tilde{A}^{\text{-}1} 
        \begin{bmatrix}
        \bm{u} \\
        \bm{0}
    \end{bmatrix}
    =
       \left[  \bm{w}^T \hspace{0.2cm} \bm{v}^T \right]     
        \begin{bmatrix}
        \bm{u}\\
        \bm{0}
    \end{bmatrix}
    = \bm{w}^T  \bm{u},
\end{equation}

 \noindent where $a_{i} = \mathcal{L}\phi \left(\|\bm{x} - \bm{x}_i\|_2 \right)|_{\bm{x}_e}$ corresponds to the entries of $\bm{a} \in \mathbb{R}^{n}$ whereas $b_{j} = \mathcal{L} p_j(\bm{x})|_{\bm{x}_e}$ corresponds to the entries of $\bm{b} \in \mathbb{R}^{\ell}$. The weights, $\bm{w} \in \mathbb{R}^{n}$, can equivalently be computed by solving the linear system,

\begin{equation} 
    \begin{bmatrix}
        A & P\\
        P^T & 0
    \end{bmatrix}
    \begin{bmatrix}
        \bm{w} \\
        \bm{v}
    \end{bmatrix}
    =
    \begin{bmatrix}
        \bm{a} \\
        \bm{b}
    \end{bmatrix},
\end{equation}

\noindent while the weights $\bm{v} \in \mathbb{R}^{\ell}$ are discarded.


\subsection{RBF-FD-C method}
\noindent Let us assume that the continuous operators appearing in (\ref{eq:peq}) have been approximated by means of RBF-FD such that $\nabla^2 \approx L$, $\mathcal{B}_{\text{D}} \approx B_{\text{D}}$ and $\mathcal{B}_{\text{N}} \approx B_{\text{N}}$. The traditional RBF-FD-C method \cite{FF2015_RBFPrimer,bayona2017role} will lead to the following square system of equations,

\begin{equation} \label{eq:eqcol}
    \begin{bmatrix}
        & & & L & & &\\
        & & & B_{\text{D}} & & &\\
        & & & B_{\text{N}} & & &
    \end{bmatrix}
    \begin{bmatrix}
        \bm{u}_{\text{I}} \\
        \bm{u}_{\text{D}}\\
        \bm{u}_{\text{N}}
    \end{bmatrix}
    =
    \begin{bmatrix}
        \bm{f} \\
        \bm{g} \\
        \bm{h} 
    \end{bmatrix},
\end{equation}

\noindent  where $L \in \mathbb{R}^{\text{N}_\text{I} \times \text{N}}$, $B_{\text{D}} \in \mathbb{R}^{\text{N}_\text{D} \times \text{N}}$, $B_{\text{N}} \in \mathbb{R}^{\text{N}_\text{N} \times \text{N}}$ and the total number of collocation nodes is $\text{N} = \text{N}_{\text{I}} + \text{N}_{\text{D}} + \text{N}_{\text{N}}$. All three discrete operators are here approximated by stencil weights based on both interior and boundary nodes according to figure \ref{fig:fitnodes}. The solution in the interior corresponds to $\bm{u}_\text{I} \in \mathbb{R}^{\text{N}_\text{I}}$, whereas the solution on the Dirichlet boundary corresponds $\bm{u}_\text{D} \in \mathbb{R}^{\text{N}_\text{D}}$ and finally, the solution on the Neumann boundary corresponds to $\bm{u}_{\text{N}} \in \mathbb{R}^{\text{N}_\text{N}}$.

\subsection{Lagrange multiplier-based RBF-FD methods}
\begin{figure}[h]
	\centering
		\includegraphics[scale=0.75]{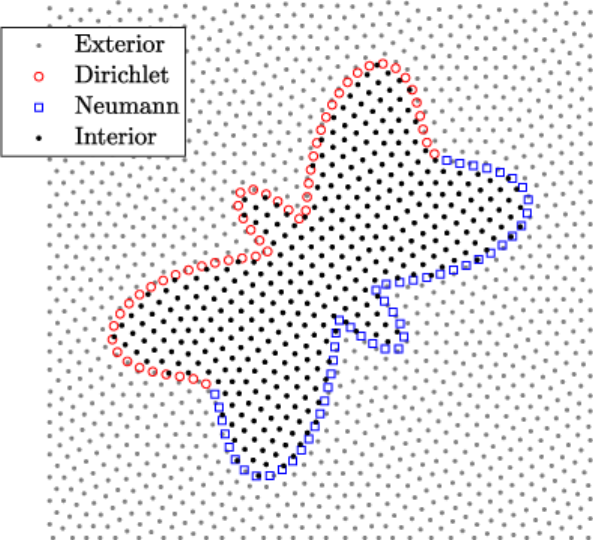} \vspace{-0cm}
	\caption{Unfitted node set including definitions of exterior (gray dots), Dirichlet (red circles), Neumann (blue squares) and interior (black dots) boundary nodes. The interior nodes are distributed with no respect to the boundary nodes.}
	\label{fig:unfitnodes}
\end{figure}

In this section, we propose for the problem in (\ref{eq:peq}) two novel optimization inspired strategies. For both, the Poisson equation is posed as a minimization problem, while the boundary conditions are enforced as constraints using the method of Lagrange multipliers. The reason for investigating these new formulations is that only nodes internal to the domain are used for computing the discrete operators, $L$, $B_{\text{D}}$ and $B_{\text{N}}$, which provides high-order approximations also when using unfitted node sets as the one illustrated in figure \ref{fig:unfitnodes}.
The two different formulations investigated are RBF-FD-LM1 and RBF-FD-LM2. These lead to two different systems of linear equations and may therefore perform differently. Both strategies are outlined below before their performance is assessed.

\subsubsection{Lagrange multiplier method for equality-constrained linear systems}

This section provides some general background that will be utilized in Sections \ref{sec:lambda1} and \ref{sec:lambda2}.

      The novelty of this work is that we pose the elliptic Poisson problem as an equality-constrained minimization problem, i.e. a norm of the residual ($A \bm{x} - \bm{b}$) is minimized while boundary conditions are enforced as equality constraints.
      Let the equality-constrained minimization problem be defined as,

\begin{equation} \label{eq:l2min}
\begin{aligned}
\min_{\bm{x}} \quad & J(\bm{x})\\
\textrm{s.t.} \quad & C \bm{x} = \bm{c}.
\end{aligned}
\end{equation}

\begin{theorem} \label{theo:t1}
The vector $\bm{x}$ that for (\ref{eq:l2min}) minimizes $J(\bm{x}) = \frac{1}{2} || A\bm{x} - \bm{b}||^2_2 = \frac{1}{2}(A\bm{x}-\bm{b})^T (A\bm{x}-\bm{b})$ can be obtained by solving the square linear system,

\begin{equation} \label{eq:con1}
    \begin{bmatrix}
        A^{T}A & C^T\\
        C & O \\
    \end{bmatrix}
    \begin{bmatrix}
        \bm{x}\\
        \bm{\lambda}\\
        \end{bmatrix}
    =
    \begin{bmatrix}
        A^{T}\bm{b} \\
        \bm{c} \\
    \end{bmatrix}.
\end{equation}

\end{theorem}

\begin{theorem} \label{theo:t2}
    If $A$ is square, symmetric, and negative definite, then the vector $\bm{x}$ that for (\ref{eq:l2min}) minimizes $J(\bm{x}) = \frac{1}{2}|| (-A)^{-1/2} ( A\bm{x} - \bm{b})||^2_2 = \frac{1}{2}(A\bm{x}-\bm{b})^T (-A)^{-1} (A\bm{x}-\bm{b})$ can be obtained by solving the square linear system,

\begin{equation} \label{eq:con2}
    \begin{bmatrix}
        A & C^T\\
        C & O \\
    \end{bmatrix}
    \begin{bmatrix}
        \bm{x}\\
        \bm{\lambda}\\
        \end{bmatrix}
    =
    \begin{bmatrix}
        \bm{b} \\
        \bm{c} \\
    \end{bmatrix}.
\end{equation}
\end{theorem}

In equations (\ref{eq:con1}) and (\ref{eq:con2}), the part $\bm{\lambda}$ of the solution vectors contain Lagrange multipliers and are typically discarded. The two novel RBF-FD-Lagrange multiplier (RBF-FD-LM) formulations (in sections \ref{sec:lambda1} and \ref{sec:lambda2}) directly use the results of Theorems \ref{theo:t1} and \ref{theo:t2}, respectively.

\subsubsection{RBF-FD-LM - Formulation 1} \label{sec:lambda1}

In the first formulation, the continuous Laplace operator is discretized at all interior nodes with stencils using only nearby interior nodes for establishing stencil weights, and then the sum of squared interior residuals are minimized. The boundary conditions are enforced as constraints on approximations that have utilized only interior nodes,

\begin{equation} \label{eq:opt1}
\begin{aligned}
\min_{\bm{u}_{\text{I}}} \quad & \frac{1}{2} || L\bm{u}_{\text{I}} - \bm{f}||^2_2\\
\textrm{s.t.} \quad & B_{\text{D}} \bm{u}_{\text{I}} = \bm{g}, \\
\quad & B_{\text{N}} \bm{u}_{\text{I}} = \bm{h}, \\
\end{aligned}
\end{equation}

\noindent where the Laplace operator ($L \in \mathbb{R}^{\text{N}_\text{I} \times \text{N}_\text{I}}$), Dirichlet boundary conditions ($B_{\text{D}} \in \mathbb{R}^{\text{N}_\text{D} \times \text{N}_\text{I}}$) and Neumann boundary conditions ($B_{\text{N}} \in \mathbb{R}^{\text{N}_\text{N} \times \text{N}_\text{I}}$) are approximated only from the solution on nodes that are considered to be interior nodes, i.e., $\bm{u}_{\text{I}} \in \mathbb{R}^{\text{N}_\text{I}}$. By Theorem \ref{theo:t1}, this can be rearranged in matrix-vector format,

\begin{equation} \label{eq:eqlam1}
    \begin{bmatrix}
        L^{T}L & B_{\text{D}}^{T} & B_{\text{N}}^{T}\\
        B_{\text{D}} & O_{\text{D,D}} & O_\text{D,N} \\
        B_{\text{N}} & O_{\text{N,D}} & O_{\text{N,N}}
    \end{bmatrix}
    \begin{bmatrix}
        \bm{u}_{\text{I}} \\
        \bm{\lambda}_{\text{D}}\\
        \bm{\lambda}_{\text{N}}
    \end{bmatrix}
    =
    \begin{bmatrix}
        L^{T}\bm{f} \\
        \bm{g} \\
        \bm{h} 
    \end{bmatrix},
\end{equation}

\noindent where $O_{\text{D,N}}$ is a zero matrix of size $\text{N}_{\text{D}} \times \text{N}_{\text{N}}$ (and similarly for the other zero matrices). Again, it must be emphasized that the discrete operators, i.e. $L$, $B_{\text{D}}$ and $B_{\text{N}}$, are here collocated only on interior nodes as depicted with black nodes in figure \ref{fig:fitnodes} and \ref{fig:unfitnodes}. Thus, both the Dirichlet and Neumann boundary conditions are essentially approximated by RBF-FD extrapolation. However, as will be illustrated in section \ref{sec:numexp}, it is still possible to use boundary-fitted node sets with/without the boundary nodes included as interior nodes. 

\subsubsection{RBF-FD-LM - Formulation 2} \label{sec:lambda2}

 In the second formulation, the Laplace operator is augmented by the boundary constraints, which mimicks (\ref{eq:con2}) from Theorem \ref{theo:t2}. The applicability of (\ref{eq:con2}) in the context it will be used here (with $L$ not necessarily symmetric, negative definite) is justified by $L$ being the discrete approximation to the continuous Laplace operator which possesses these properties.
The formulation gives rise to the system of equations,

\begin{equation} \label{eq:eqlam2}
    \begin{bmatrix}
        L & B_{\text{D}}^{T} & B_{\text{N}}^{T}\\
        B_{\text{D}} & O_{\text{D,D}} & O_\text{D,N} \\
        B_{\text{N}} & O_{\text{N,D}} & O_{\text{N,N}}
    \end{bmatrix}
    \begin{bmatrix}
        \bm{u}_{\text{I}} \\
        \bm{\lambda}_{\text{D}}\\
        \bm{\lambda}_{\text{N}}
    \end{bmatrix}
    =
    \begin{bmatrix}
        \bm{f} \\
        \bm{g} \\
        \bm{h} 
    \end{bmatrix},
\end{equation}

\noindent where all submatrices are the equivalent to the ones used in  (\ref{eq:eqlam1}), however, the matrix-matrix product $L^{T}L$ and matrix-vector product $L^{T} \bm{f}$ have been replaced by $L$ and $\bm{f}$, respectively. Other references to using Lagrange multipliers as constraints to the Poisson problem can be found in e.g. \cite{trottenberg2000multigrid,wright2023mgm}.

\subsubsection{Comments on the RBF-FD-LM methods}
\noindent Different possibilities arise when the problem in (\ref{eq:peq}) is formulated as in (\ref{eq:eqlam1}) or (\ref{eq:eqlam2}). Firstly, it is possible to use the same boundary-fitted node set as used in the traditional RBF-FD-C method as illustrated in figure \ref{fig:fitnodes}. Secondly, it is also possible to use an unfitted node set as illustrated in figure \ref{fig:unfitnodes}, and use only the interior nodes as collocation nodes. Thus, the discrete boundary conditions will be approximated at $\text{N}_{\text{D}} + \text{N}_{\text{N}}$ nodes on $\Gamma$ from a subset of the interior nodes using the exact same algorithm for computing RBF-FD weights as in the RBF-FD-C method. Hence, the only change is how the problem is posed, and the problem size $\text{N}$ stays the same as in the RBF-FD-C method if an unfitted node set is used.\\


\section{Numerical experiments} \label{sec:numexp}
These test problems are introduced to investigate the accuracy and robustness of the two novel RBF-FD-LM methods. For all test problems the relative $\ell_2$-norm error is used and defined as,

\begin{equation}
e = \frac{\parallel u_{num} - u_{exact} \parallel_2}{\parallel u_{exact} \parallel_2},
\end{equation}

\noindent where $u_{num}$ is the numerical solution and $u_{exact}$ is the exact analytical solution. In this work, the emphasis is put on showing the basic principles and only simplistic tests based on the method of manufactured solutions will be considered for now.

\subsection{Test problem 1: Poisson's equation in the unit disk with Dirichlet boundary conditions}

In this test problem, the relation between stencil sizes and augmented polynomials is investigated in a similar manner as in \cite{bayona2017role}. The test problem is the same as in \cite{bayona2017role}, and is posed as

\begin{equation}\label{eq:test1}
\begin{aligned}
\nabla^2 u &= f(x,y),\hspace{0.25cm} \bm{x} \in \Omega, \\
u  &= g(x,y), \hspace{0.25cm} \bm{x} \in \Gamma_{\text{D}},\\
\end{aligned}
\end{equation}

\noindent where $f(x,y)$ and $g(x,y)$ are computed from the exact solution $u(x,y) = \text{sin}(10(x+y))$, while the domain is defined as the unit disk $\Omega = \{ (x,y) \in \mathbb{R}^2 : x^2 + y^2 \leq 1 \}$ and the Dirichlet boundary as $\Gamma_{\text{D}} = \{(x, y) \in \mathbb{R}^2 : {x^2 + y^2 = 1} \}$. 

\begin{figure}[h!]
	\centering
		\includegraphics[scale=0.4]{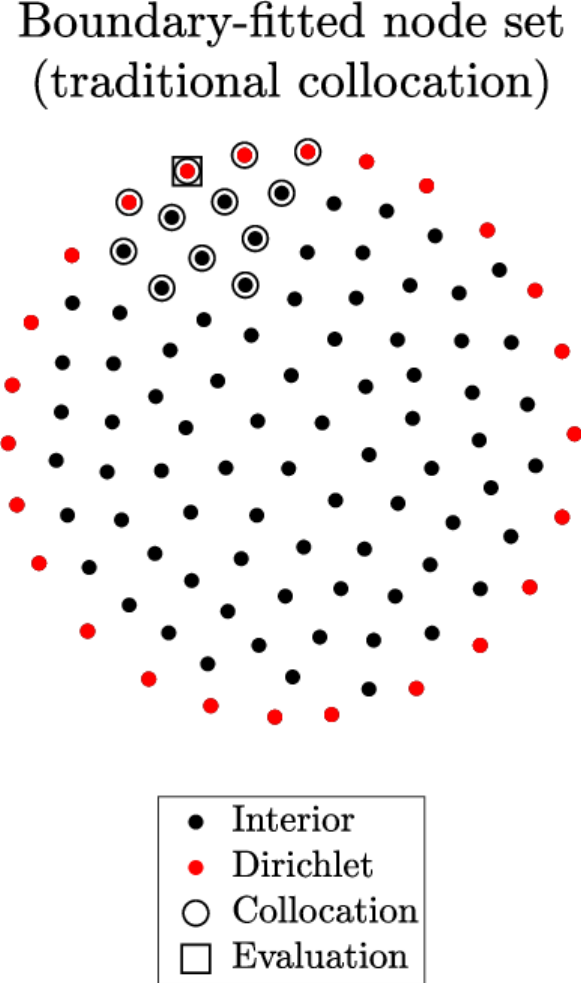} 
  		\hspace{0.2cm} 
    \includegraphics[scale=0.4]{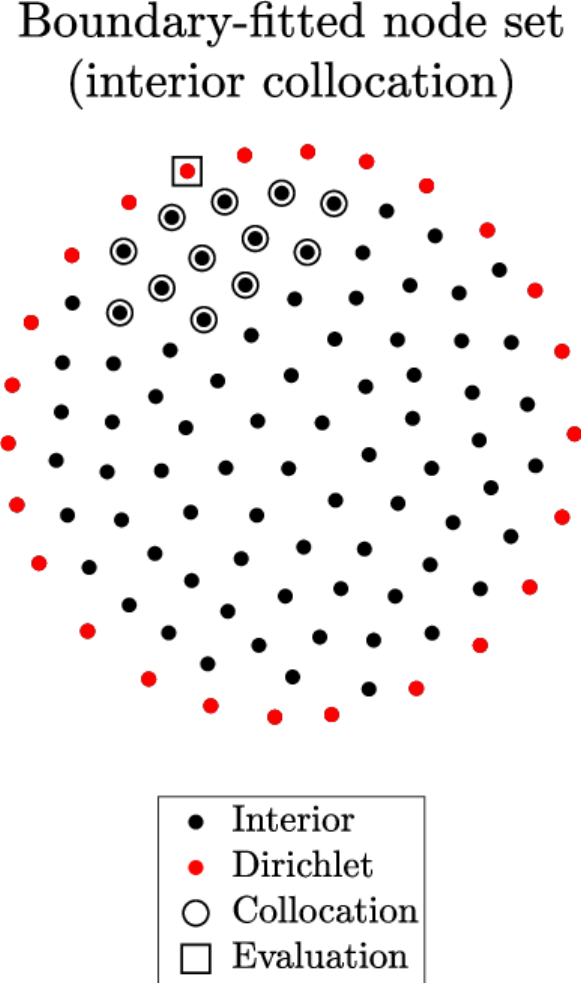} \hspace{0.2cm}
  \includegraphics[scale=0.4]{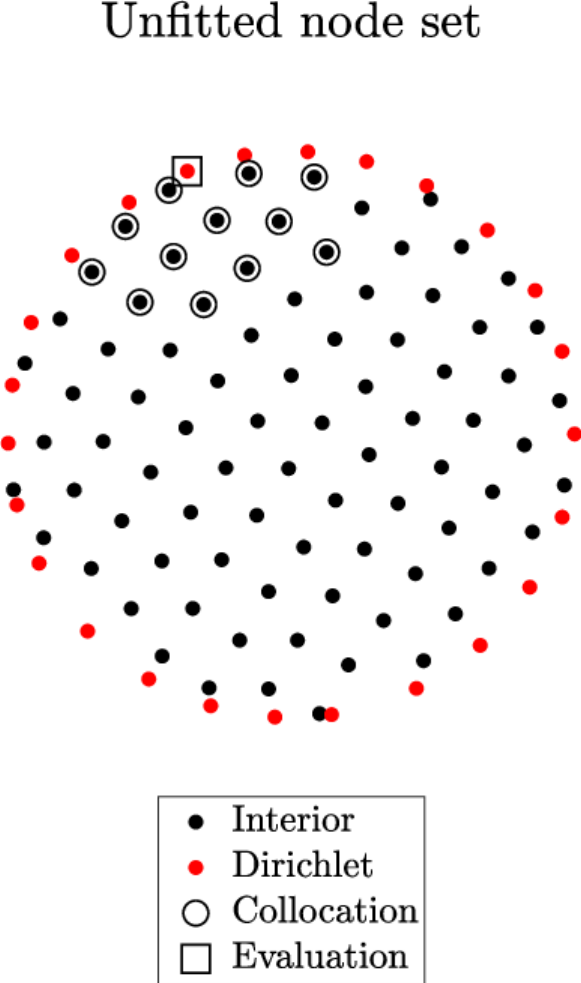}
	\caption{The three different node sets and collocation setups used for test problem 1. The stencil illustrated on each node set is based on a relative stencil size $\lceil n/\ell \rceil = 2$ and augmented polynomial degree $m = 2$.}
	\label{fig:nset_test1}
\end{figure}

The Poisson problem in  (\ref{eq:test1}) is solved using the traditional RBF-FD-C method and the new RBF-FD-LM methods. This initial test seeks to illustrate the consistency and robustness of the RBF-FD-LM methods. The test problem is solved on three different node sets and collocation setups, as illustrated in figure \ref{fig:nset_test1}, which serves to illustrate the versatility of the RBF-FD-LM methods. In all three node set configurations, the interior nodes in $\Omega$ are determined from a node set generated in a square domain covering $\Omega$ using the algorithm proposed in \cite{fornberg2015fast}, while a node repelling algorithm similar to the one in \cite{bayona2017role} has been used to establish the boundary-fitted node sets.

\begin{figure}[h!]
	\centering
		\includegraphics[scale=0.6]{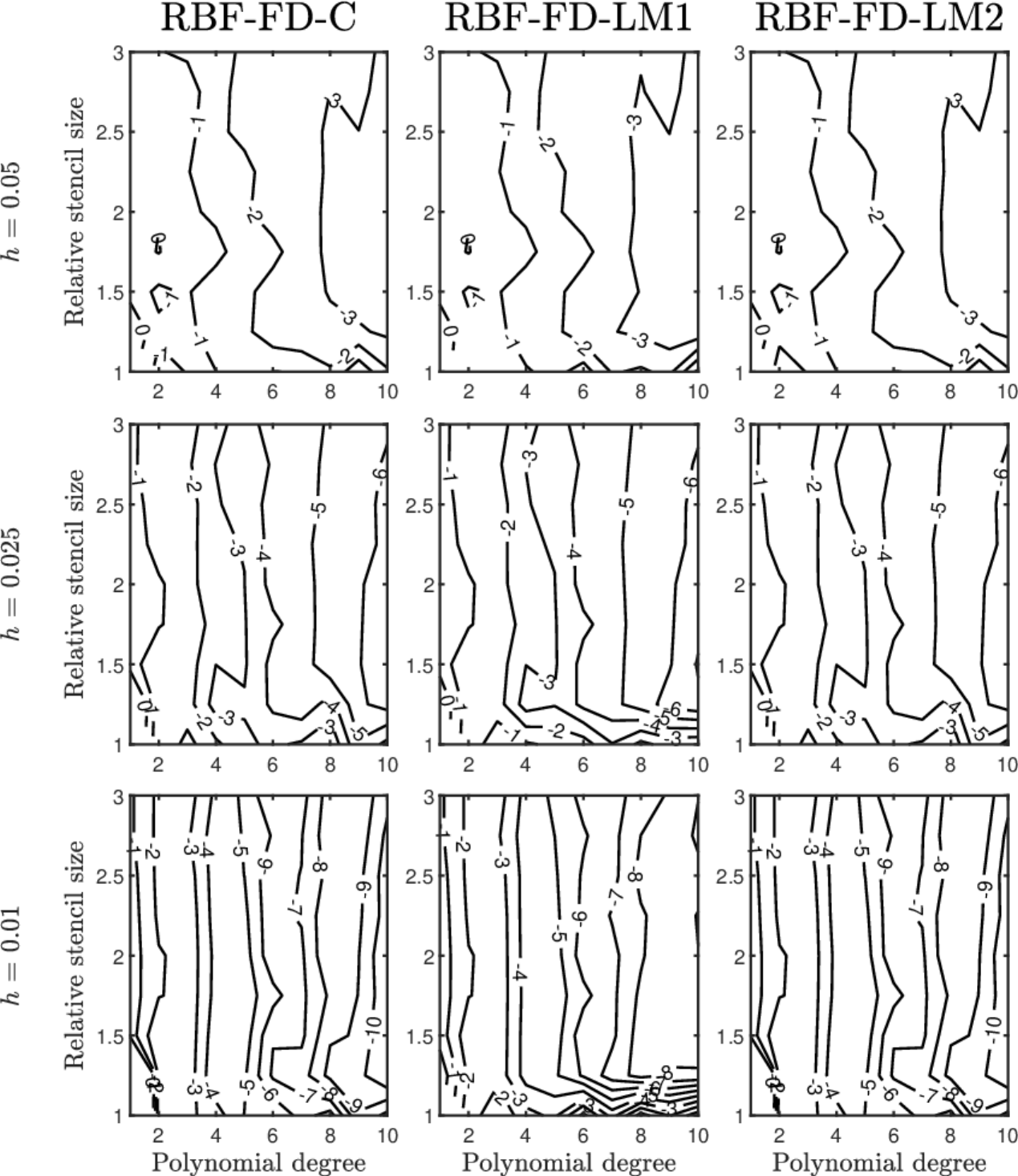}
	\caption{$\text{Log}_{10}$ of the relative $\ell_2$-norm errors as function of relative stencil size ($\lceil n/\ell \rceil$) and augmented polynomial degree ($m$) using the radial basis function $\phi(r) = r^3$ on a boundary-fitted node set, where all three methods use interior and boundary nodes as collocation nodes, for solving the Poisson equation in the unit disk with Dirichlet boundary conditions.}
	\label{fig:cont_fit}
\end{figure}

In figure \ref{fig:cont_fit}, the relative $\ell_2$-norm errors are shown for the three different RBF-FD methods as function of augmented polynomial degree ($m$) and relative stencil size ($\lceil n/\ell \rceil$), i.e. ratio between neighbors and number of terms in the augmented polynomial, while using the radial basis function $\phi(r) = r^3$ on a node set corresponding to the boundary-fitted node set with traditional collocation in figure \ref{fig:nset_test1}.

Figure \ref{fig:cont_fit} differs from its counterpart figure 2 in \cite{bayona2017role} as the vertical axes show relative stencil size ($n/ \ell$) rather than stencil size ($n$). As background to considering $n/ \ell$ , we note that, in $d$ dimensions, there are $\ell = \binom{d+m}{m}$ independent polynomials up through degree $m$. As was observed in \cite{bayona2017role}, the system for obtaining the RBF-FD 
weights is singular if $n/ \ell$ $<$ 1 and reliable accuracy requires $n/ \ell$ $\geq$ 2. Hence, we choose here to use $n/ \ell$ on the 
vertical axis (rather than $n$) and focus our attention on the range 2 $\leq n/ \ell \leq$ 3. It was also noted previously that the 
polynomial degree $m$ alone controls the order of convergence under node refinement. Instead of $m$ = 7 (in Figure 2 
in \cite{bayona2017role}), we use here $m$ = 3.


Furthermore, it can be noticed that RBF-FD-LM2 method seems to obtain very similar results to the traditional RBF-FD-C method. Although, the RBF-FD-LM1 method seems to provide similar accuracy as the other two methods, it can be noticed that as the node density $h$ is decreased the solution stagnates around $e \approx \mathcal{O}(10^{-9})$, whereas the two other methods do not.

\begin{figure}[h!]
	\centering

		\includegraphics[scale=0.6]{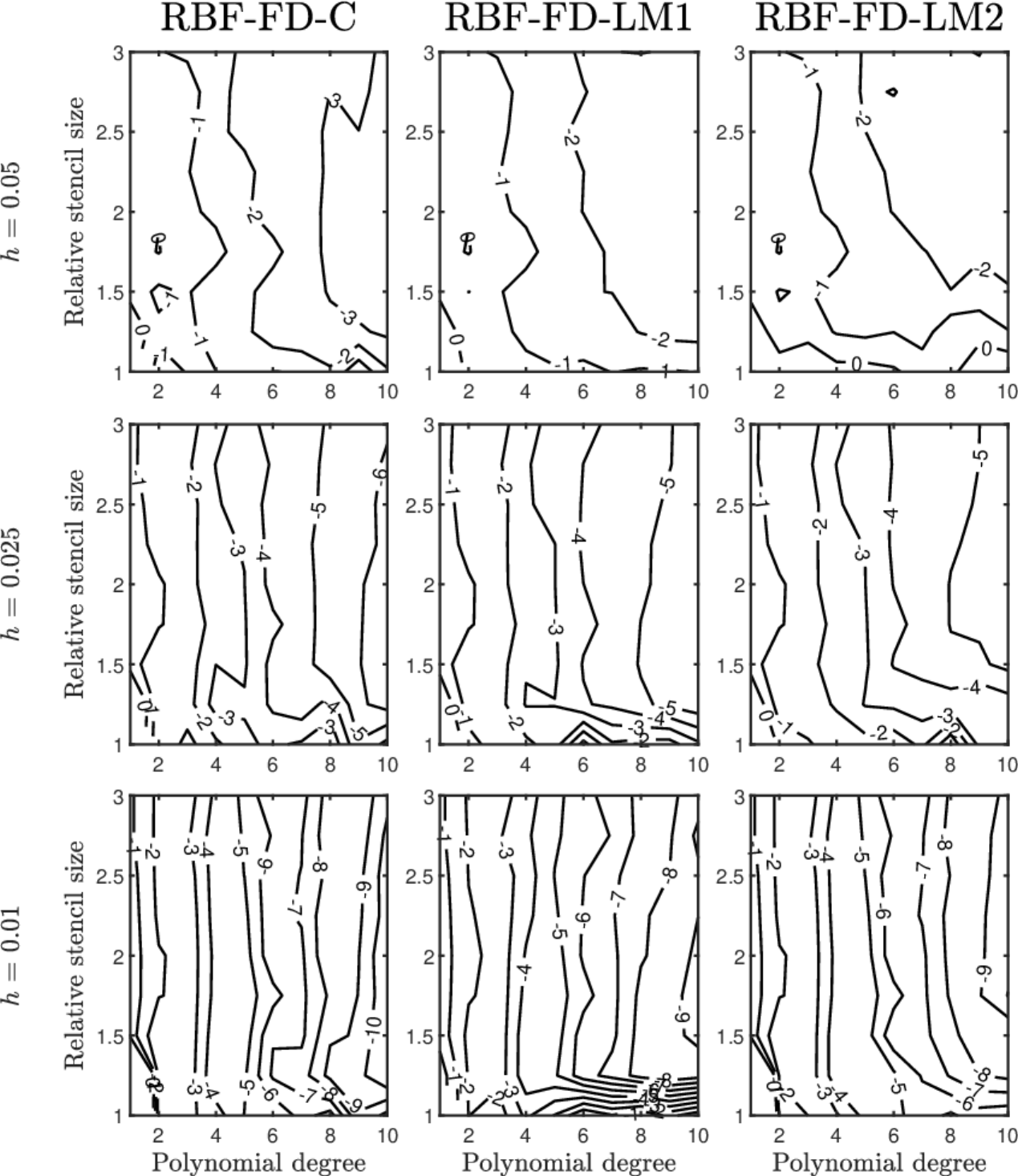}
	\caption{$\text{Log}_{10}$ of the relative $\ell_2$-norm errors as function of relative stencil size ($\lceil n/\ell \rceil$) and augmented polynomial degree ($m$) using the radial basis function $\phi(r) = r^3$ on a boundary-fitted node set, where the RBF-FD-LM methods use only interior nodes as collocation nodes, for solving the Poisson equation in the unit disk with Dirichlet boundary conditions.}
	\label{fig:cont_fit_in}
\end{figure}

Next, the problem in (\ref{eq:test1}) is again solved, but now only the interior nodes are used as collocation nodes for the two RBF-FD-LM methods, whereas the RBF-FD-C method still uses both interior and boundary nodes as collocation nodes. The errors are illustrated in figure \ref{fig:cont_fit_in}. Here, it can be noticed that the accuracy decreases about one order of magnitude for the two RBF-FD-LM methods, while the relation between the relative stencil size and augmented polynomial degree seems to follow the same trend as the traditional collocation setup.

\begin{figure}[h!]
	\centering
		\includegraphics[scale=0.5]{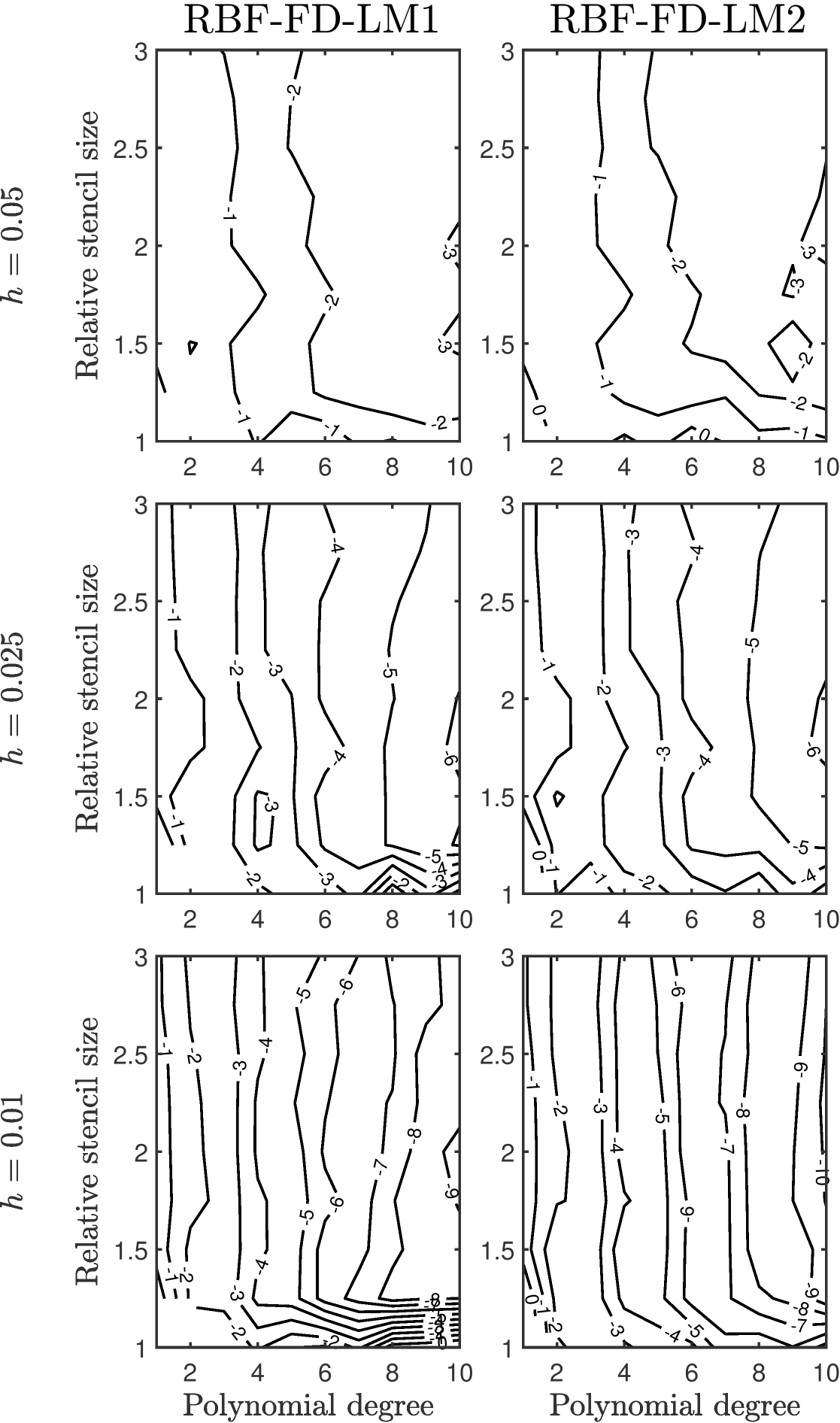}
	\caption{$\text{Log}_{10}$ of the relative $\ell_2$-norm errors as function of relative stencil size ($\lceil n/\ell \rceil$) and augmented polynomial degree ($m$) using the radial basis function $\phi(r) = r^3$ on an unfitted node set, where only interior nodes are used as collocation nodes, for solving the Poisson equation in the unit disk with Dirichlet boundary conditions.}
	\label{fig:cont_unfit}
\end{figure}

Finally, the unfitted node set is tested and as it is only applicable to the two RBF-FD-LM methods, the relative error plots are illustrated only for those two methods on figure \ref{fig:cont_unfit}.
As the interior nodes are now placed closer to the boundary than in the boundary-fitted node set with interior collocation, only a minor decrease of accuracy is noticed when compared to the traditional collocation method on a boundary-fitted node set.

\clearpage
\newpage

\subsection{Test problem 2: Poisson's equation in butterfly domain with Dirichlet and Neumann boundary conditions}

\begin{figure}[h!]
	\centering
         \includegraphics[scale=0.625]{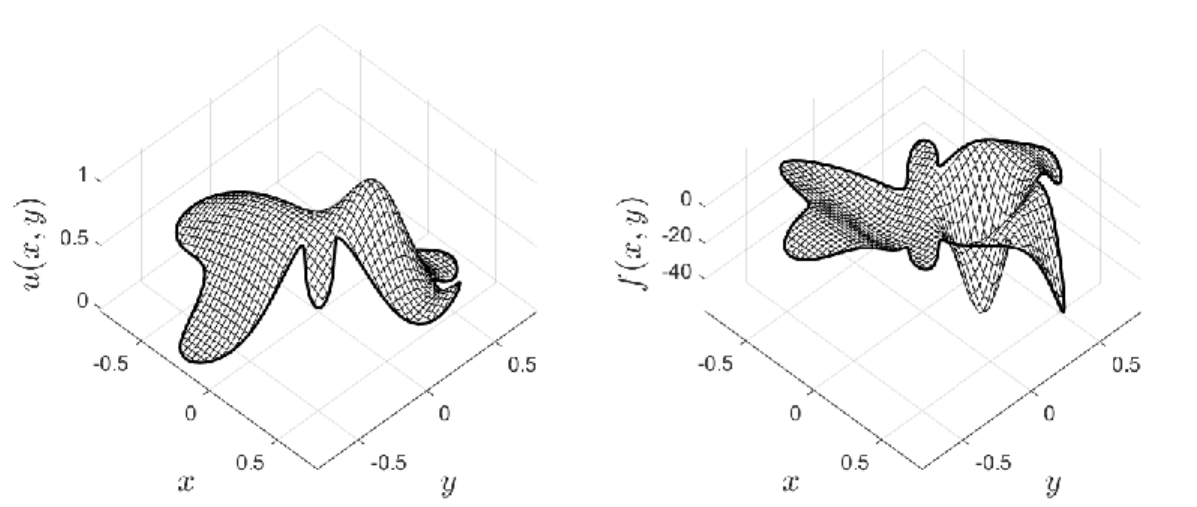}
	\caption{The exact solution and source function for test problem 2.}
	\label{fig:unfit_DN_u_f}
\end{figure}

This test considers a more complicated two-dimensional domain with both Dirichlet and Neumann boundary conditions, which was also used in \cite{tominec2021unfitted}. The test problem is expressed as,

\begin{equation}\label{eq:test2}
\begin{aligned}
\nabla^2 u &= f(x,y),\hspace{0.25cm} \bm{x} \in \Omega, \\
u  &= g(x,y), \hspace{0.25cm} \bm{x} \in \Gamma_{\text{D}},\\
\frac{\partial u}{\partial n}  &= h(x,y), \hspace{0.25cm} \bm{x} \in \Gamma_{\text{N}},\\
\end{aligned}
\end{equation}

\noindent and the domain $\Omega$ is the butterfly illustrated in figure \ref{fig:unfitnodes}, which is represented in polar coordinates as,

\begin{equation}
r(\theta) = \frac{1}{4} \left( 2 + \text{sin}(2\theta) - 0.01\text{cos} \left(5 \theta - \frac{\pi}{2}\right) + 0.63 \text{sin}(6 \theta - 0.1)\right),
\end{equation}

\noindent where $f(x,y)$, $g(x,y)$ and $h(x,y)$ are computed from the exact solution, 

\begin{equation}
\begin{aligned}
u(x,y) =&  \frac{3}{4}e^{-\frac{1}{4}((9x-2)^2 + (9y-2)^2)} + \frac{3}{4}e^{-(\frac{1}{49}(9x+1)^2 + \frac{1}{10}(9y+1)^2)}\\
    &+ \frac{1}{2}e^{-\frac{1}{4}((9x-7)^2 + (9y-3)^2)} - \frac{1}{5}e^{-((9x-4)^2 + (9y-7)^2)},
\end{aligned}
\end{equation}

\noindent also known as the Franke function. The interior nodes in $\Omega$ are determined from a node set generated in a square domain covering $\Omega$ using the algorithm proposed in \cite{fornberg2015fast}. An example of the node distribution used is illustrated in figure \ref{fig:unfitnodes}, while the non-trivial source function is illustrated in figure \ref{fig:unfit_DN_u_f} along with the exact solution.

\begin{figure}[h!]
	\centering
		\includegraphics[scale=0.35]{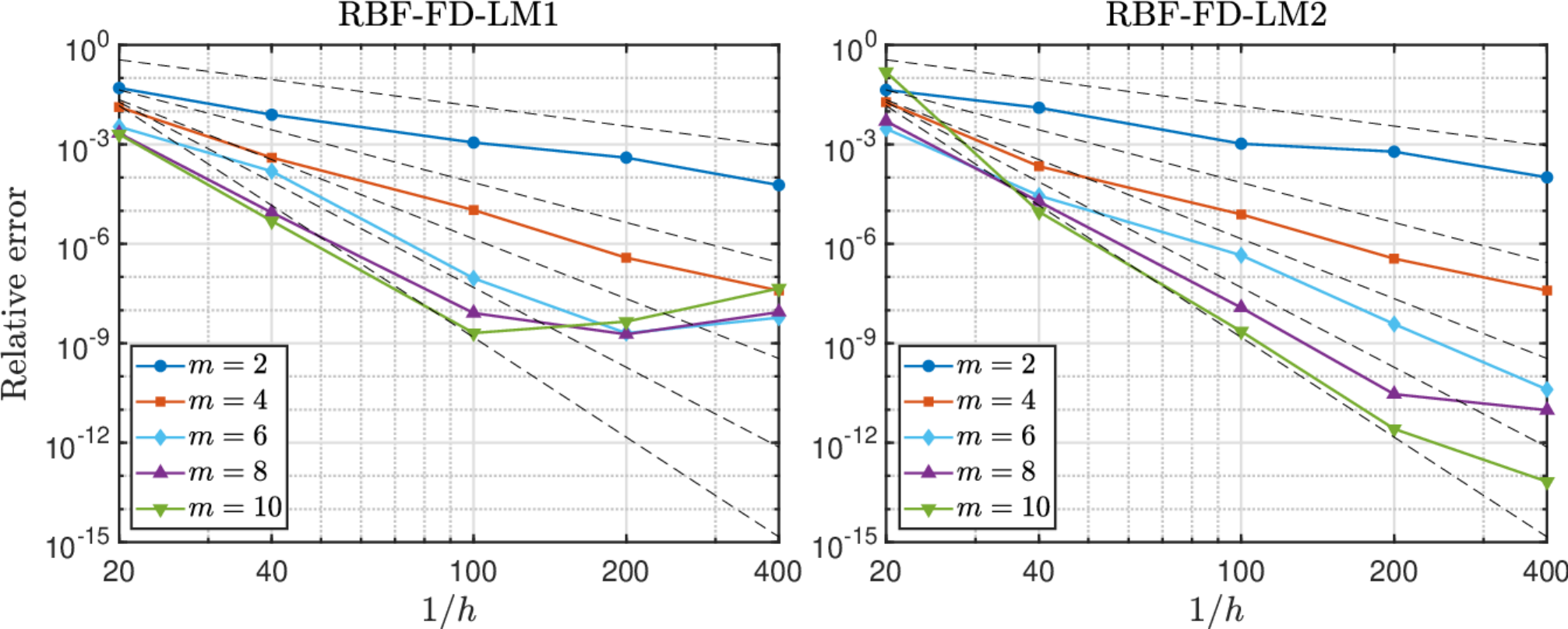}
	\caption{Convergence rates for both RBF-FD-LM methods with $\phi(r) = r^3$, varying augmented polynomial degrees ($m$) and relative stencil sizes $\lceil n/\ell \rceil = 2$. The dashed lines indicate slopes equal to 2, 4, 6, 8 and 10, respectively.}
	\label{fig:conv_unfit_DN}
\end{figure}

The test problem in (\ref{eq:test2}) is solved using $\phi(r) = r^3$ and varying polynomial degrees ($m$) for both RBF-FD-LM methods. Convergence rates are illustrated in figure \ref{fig:conv_unfit_DN}.

From figure \ref{fig:conv_unfit_DN} it can be seen that the convergence rates of both  RBF-FD-LM methods are similar until formulation 1 stagnates around $e \approx \mathcal{O}(10^{-9})$, which was also noticed in test problem 1, whereas formulation 2 continuous until machine precision. The stagnation of formulation 1 is is not altogether unexpected, as this case contains a matrix block $L^{T}L$ in place of $L$, which for square matrices squares the condition number.

\begin{figure}[h!]
	\centering
		\includegraphics[scale=0.35]{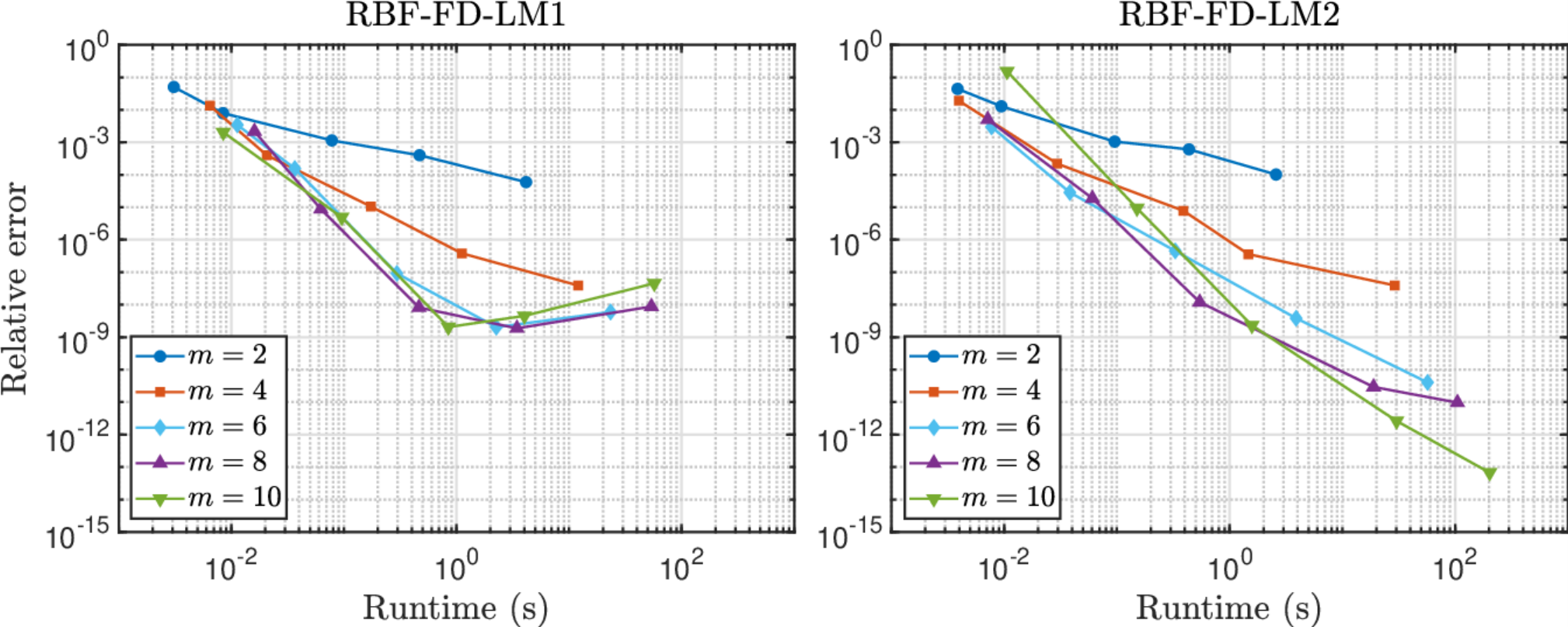}
	\caption{Computational efficiency for both RBF-FD-LM methods with $\phi(r) = r^3$, varying augmented polynomial degrees ($m$) and relative stencil sizes $\lceil n/\ell \rceil = 2$.}
	\label{fig:conv_unfit_DN_cpu}
\end{figure}

In addition to the convergence rates, the computational efficiency of each RBF-FD-LM method is illustrated in figure \ref{fig:conv_unfit_DN_cpu}. It must be noted that all computations are carried out in MATLAB with \textit{mldivide()} and the runtime includes, in addition to the solve, the time to assemble the matrices in  (\ref{eq:eqlam1}) and (\ref{eq:eqlam2}), while it excludes the RBF-FD weight computations. The computational efficiency is of major importance for applying these methods to large-scale problems, however, we leave out these investigations for the sake of conciseness.

\subsection{Test problem 3: Poisson's equation in the sphere with Dirichlet boundary conditions}

\begin{figure}[h!]
	\centering
		\includegraphics[scale=0.35]{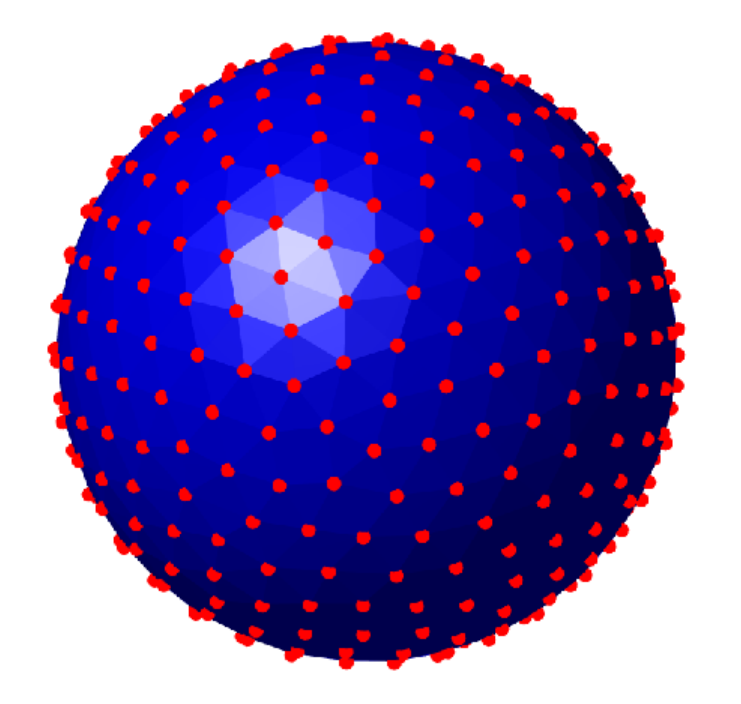}
	\caption{A visualization of the Dirichlet boundary nodes used for test problem 3 (coarsest node set level in the convergence study).}
	\label{fig:nodes_3d_ex}
\end{figure}

The third test problem considers a simple three-dimensional domain, namely the unit sphere, with Dirichlet boundary conditions. The test problem is expressed as,

\begin{equation}\label{eq:test3}
\begin{aligned}
\nabla^2 u &= f(x,y,z),\hspace{0.25cm} \bm{x} \in \Omega, \\
u  &= g(x,y,z), \hspace{0.25cm} \bm{x} \in \Gamma_{\text{D}},\\
\end{aligned}
\end{equation}

\noindent where $f(x,y,z)$ and $g(x,y,z)$ are computed from the exact solution $u(x,y,z) = \text{sin}(\pi x) + \text{cos}(\pi y) + \text{sin}(\pi z)$, while the domain is defined as the closed unit ball $\Omega = \{ (x,y,z)\in \mathbb{R}^3, x^2 + y^2 + z^2 \leq 1 \}$. The interior nodes in $\Omega$ have been generated using the algorithm presented in \cite{van2021fast} to scatter nodes in a cube that covers the domain $\Omega$ in three-dimensional space. For the sake of simplicity, we have used pre-computed surface node sets to represent the boundary $\Gamma_{\text{D}} = \{(x, y, z) \in \mathbb{R}^3 : {x^2 + y^2 + z^2 = 1} \}$, which is desribed in \cite{sloan2004extremal} and made available at \cite{robnodes}. An example of the boundary node distribution is illustrated in figure \ref{fig:nodes_3d_ex}.

\begin{figure}[h!]
	\centering
		\includegraphics[scale=0.35]{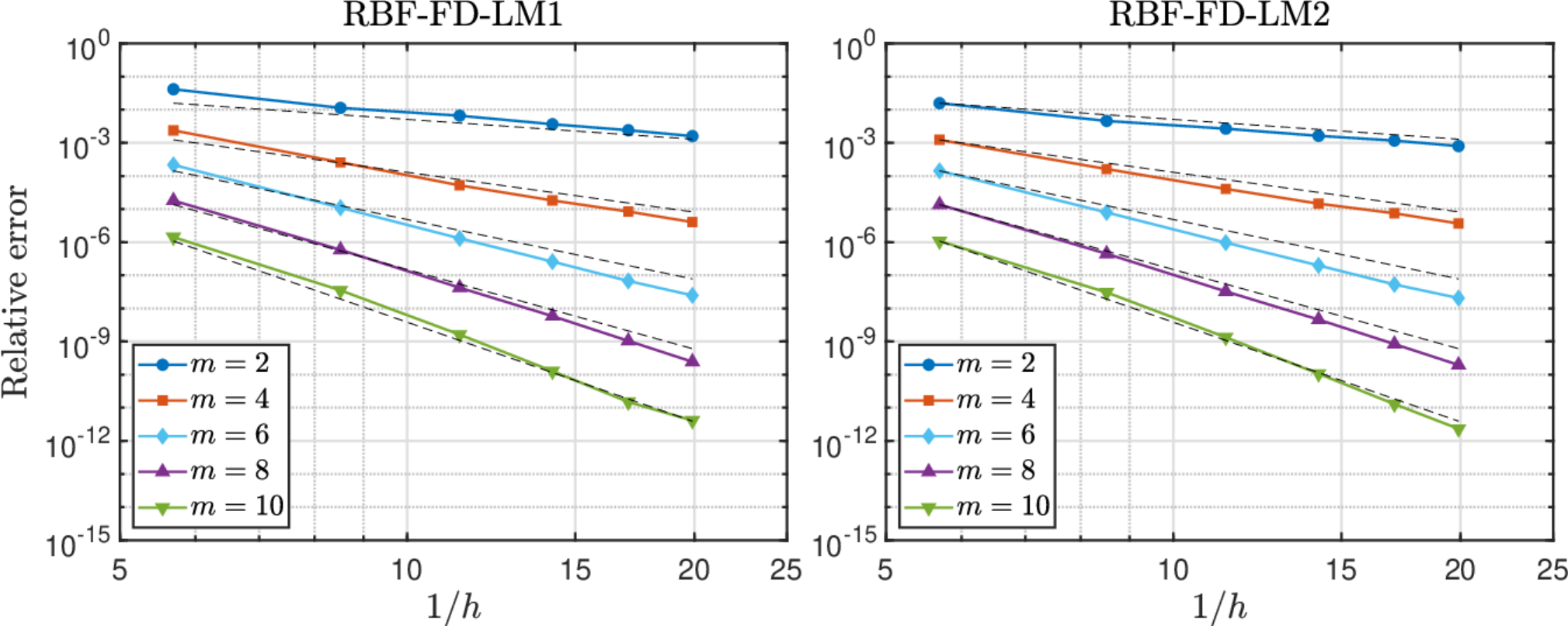}
	\caption{Convergence rates for both RBF-FD-LM methods applied to test problem 3 using $\phi(r) = r^3$, varying augmented polynomial degrees ($m$) and relative stencil sizes $\lceil n/\ell \rceil = 2$. The dashed lines indicate slopes equal to 2, 4, 6, 8 and 10, respectively.}
	\label{fig:conv_unfit_D_sphere}
\end{figure}

Again, the test problem in (\ref{eq:test3}) is solved using $\phi(r) = r^3$ and varying polynomial degrees ($m$) for both RBF-FD-LM methods and the relative stencils sizes are again $\lceil n/\ell \rceil = 2$. Convergence rates are illustrated in figure \ref{fig:conv_unfit_D_sphere}.\\

From the results in figure \ref{fig:conv_unfit_D_sphere}, it is again seen how both formulations provide high-order accurate solutions. However, formulation 1 does not seem to run into any stagnation errors, which were noticed for test problem 2, and the formulations obtain very similar orders of accuracy. It must be noted here, that the node density is not as fine as in the previous two-dimensional test problems, which could be the reason for the presence of stagnation errors noticed for formulation 1.

\subsection{Test problem 4: Poisson's equation in the sphere with Dirichlet and Neumann boundary conditions}

\begin{figure}[h!]
	\centering
		\includegraphics[scale=0.45]{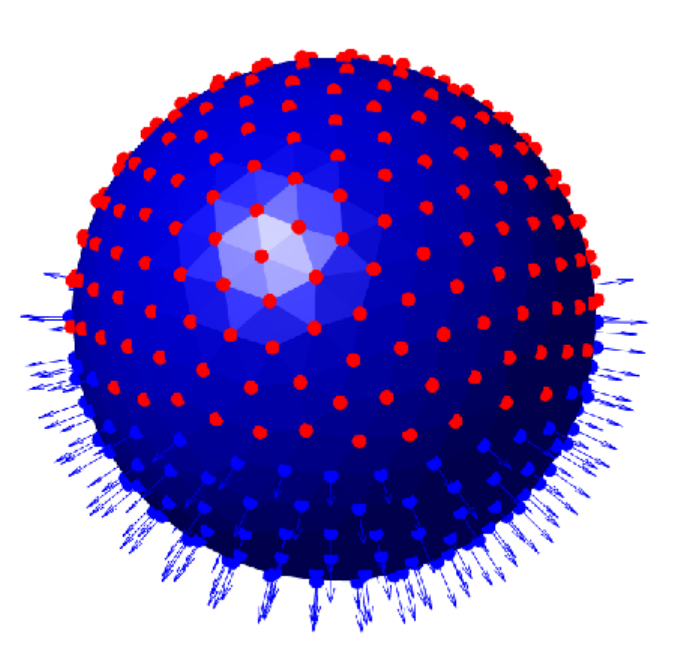}
	\caption{A visualization of the Dirichlet (red points) and Neumann (blue points and arrows) boundary nodes used for test problem 4 (showing the coarsest node set in the convergence study).}
	\label{fig:nodes_3d_ex_DN} 
\end{figure}

The last test problem considers the unit sphere with both Dirichlet and Neumann boundary conditions. The test problem is expressed as,

\begin{equation}\label{eq:test4}
\begin{aligned}
\nabla^2 u &= f(x,y,z),\hspace{0.25cm} \bm{x} \in \Omega, \\
u  &= g(x,y,z), \hspace{0.25cm} \bm{x} \in \Gamma_{\text{D}},\\
\frac{\partial u}{\partial n} &= h(x,y,z), \hspace{0.25cm} \bm{x} \in \Gamma_{\text{N}},\\
\end{aligned}
\end{equation}

\noindent where $f(x,y,z)$, $g(x,y,z)$ and $h(x,y,z)$ are computed from the exact solution $u(x,y,z) = \text{sin}(\pi x) + \text{cos}(\pi y) + \text{sin}(\pi z)$, while the domain is defined as the closed unit ball $\Omega = \{ (x,y,z) \in \mathbb{R}^3, x^2 + y^2 + z^2 \leq 1 \}$. The boundary and interior nodes have been generated with the same methodology as for test problem 3, while the boundary nodes are chosen such that $\Gamma_{\text{D}} = \{(x, y, z) \in \mathbb{R}^3 : {x^2 + y^2 + z^2 = 1}, {z\geq0} \}$ and $ \Gamma_{\text{N}} = \{(x, y, z) \in \mathbb{R}^3 : {x^2 + y^2 + z^2 = 1}, {z < 0}\}$. An example of the boundary node distribution along with node definitions is illustrated in figure \ref{fig:nodes_3d_ex_DN}.

\begin{figure}[h!]
	\centering
		\includegraphics[scale=0.35]{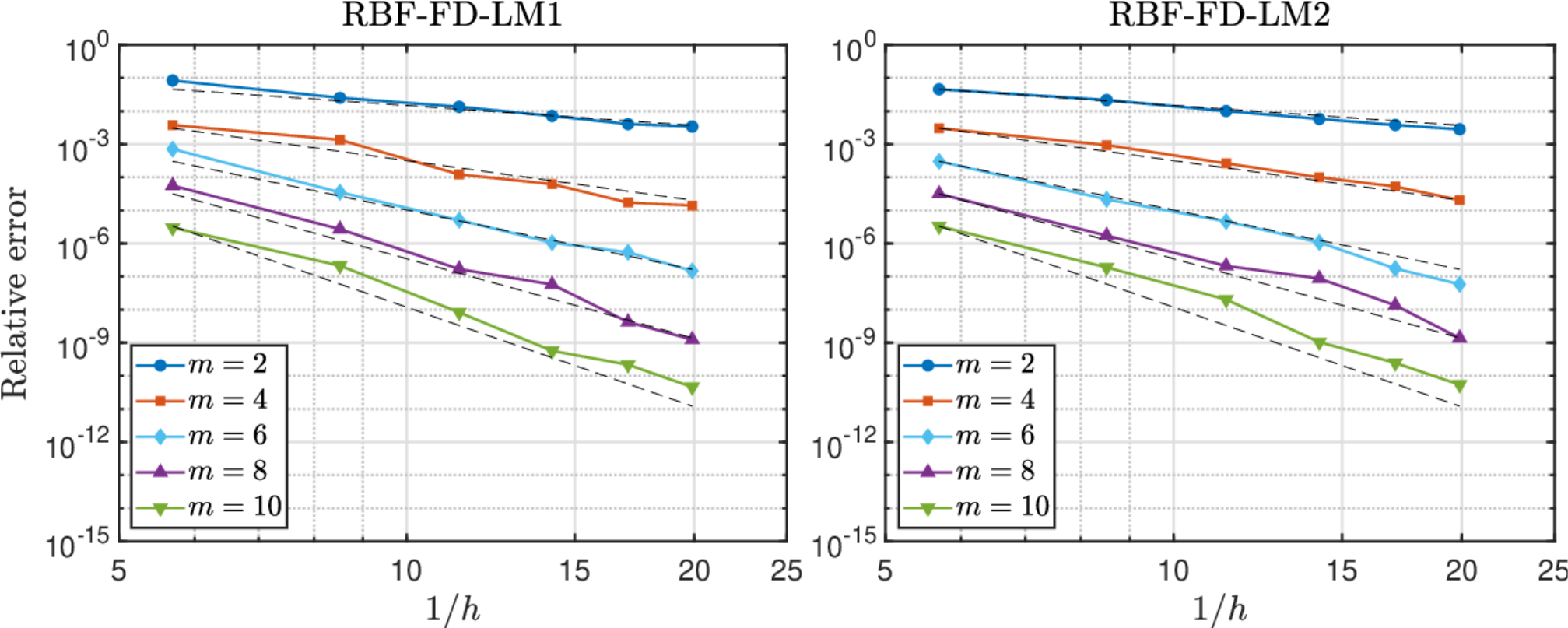}
	\caption{Convergence rates for both RBF-FD-LM methods applied to test problem 4 using $\phi(r) = r^3$, varying augmented polynomial degrees ($m$) and relative stencil sizes $\lceil n/\ell \rceil = 2$. The dashed lines indicate slopes equal to 2, 4, 6, 8 and 10, respectively.}
	\label{fig:conv_unfit_DN_sphere}
\end{figure}

The test problem in (\ref{eq:test4}) is solved using $\phi(r) = r^3$ and varying polynomial degrees ($m$) for both RBF-FD-LM methods and the relative stencil sizes are again $\lceil n/\ell \rceil = 2$. Convergence rates are illustrated in figure \ref{fig:conv_unfit_DN_sphere}.

As seen from the results in figure \ref{fig:conv_unfit_DN_sphere}, high-order accurate solutions are obtained for both RBF-FD-LM methods and no error stagnation seems to be present in either of the formulations. However, the node density is not as refined as in the two-dimensional test problems. 

Finally, it must be emphasized that the same relative stencil sizes have been used for both the two- and three-dimensional test problems. Similar results can possibly be achieved with smaller relative stencil sizes especially in three-dimensional cases.

\section{Conclusions}
\label{sec:concl}

In this paper, a novel method for solving elliptic PDEs discretized using RBF-FD has been presented, which maintains high-order accuracy on unfitted node sets. This is achieved by combining Lagrange multiplier boundary constraints with pure interior collocation of the RBF-based interpolants, which has led to the two different formulations proposed here.

Differences in the derivation of the two proposed formulations have been outlined and numerical experiments illustrate that high-order accurate solutions of Poisson's equation can be achieved for both two- and three-dimensional problems with Dirichlet and Neumann boundary conditions.

Ideas presented in this paper, and implemented with RBF-FD, can potentially also be implemented using other spatial discretization schemes, and preliminary tests have shown a potential for applications within time dependent (e.g. moving boundary) problems as well.

\bibliographystyle{siamplain}
\bibliography{main}
\end{document}